\documentclass[submission]{FPSAC2022}

\bibliography{positroid.bib}

 \usepackage{subcaption}
 
 \usepackage{MnSymbol}
 
\usepackage{tikz}
\usepackage{tikz-cd}
\usetikzlibrary{arrows.meta}
\usetikzlibrary{bending}

\theoremstyle{plain}
\newtheorem{Th}{Theorem}[section]
\newtheorem{Lemma}[Th]{Lemma}

\newtheorem{ThmDef}[Th]{Theorem/Definition}

\theoremstyle{remark}

\theoremstyle{definition}
\newtheorem{Def}[Th]{Definition}

\newtheorem{?}[Th]{Problem}
\newtheorem{Ex}[Th]{Example}

\newcommand{\C}{\mathbb{C}}

\newcommand{\M}{\mathcal{M}}

\newcommand{\QQ}{\mathcal{Q}}

\newcommand{\doubleloop}{\begin{tikzpicture}[scale=1.5,cap=round,>=latex]
  \def\Radius{.1cm}

  \draw (0cm,0cm) circle[radius=\Radius];

  \begin{scope}[
    -{Stealth[round, length=2pt, width=2pt, bend]},
    shorten >=1pt,
    very thin,
  ]
    \draw (\Radius, 0) arc(-3:3:\Radius);
    \draw (-\Radius, 0) arc(180+3:180-3:\Radius);
  \end{scope}
  
\end{tikzpicture}}

\newcommand\tinydoubleloop{\vcenter{\hbox{\scalebox{0.5}{\doubleloop}}}}

\newcommand{\w}{w^{\hspace{0.01in} \tinydoubleloop}}
\newcommand{\vv}{v^{\hspace{0.01in} \tinydoubleloop}}
\newcommand{\SSn}{S_n^{\hspace{0.01in} \tinydoubleloop}}

\newcommand{\SSnminus}{S_{n-1}^{\hspace{0.01in} \tinydoubleloop}}
\newcommand{\SSnk}{S_{n,k}^{\hspace{0.01in} \tinydoubleloop}}

\newcommand{\Gkn}{Gr(k,n)}
\newcommand{\aligns}{Alignments}
\newcommand{\clockwise}[1]{\overrightarrow{#1}}
\newcommand{\counterclockwise}[1]{\overleftarrow{#1}}

\DeclareMathOperator{\GL}{GL}

\DeclareMathOperator{\rank}{rank}

\usepackage{lipsum}

\title[Pattern Avoidance and Positroid Varieties]{A Pattern Avoidance
Characterization for Smoothness of Positroid Varieties}

\author[S. Billey and J. Weaver]{Sara
C. Billey\thanks{Both authors
were partially supported by the National Science Foundation grant
DMS-1764012. Email:billey@math.washington.edu and jeweaver@uw.edu.}\addressmark{1} \and Jordan E. Weaver\addressmark{1}}

\address{\addressmark{1}Department of Mathematics, University of Washington, Seattle, WA, USA}

\received{November 24, 2021}

\abstract{ Positroids are certain representable matroids originally
studied by Postnikov in connection with the totally nonnegative
Grassmannian and now used widely in algebraic combinatorics.  The
positroids give rise to determinantal equations defining positroid
varieties as subvarieties of the Grassmannian variety. Rietsch,
Knutson--Lam--Speyer and Pawlowski studied geometric and cohomological
properties of these varieties.  In this paper, we continue the study
of the geometric properties of positroid varieties by establishing
several equivalent conditions characterizing smooth positroid
varieties using a variation of pattern avoidance defined on decorated
permutations, which are in bijection with positroids. Furthermore, we
give a combinatorial method for determining the dimension of the
tangent space of a positroid variety at key points using an induced
subgraph of the Johnson graph.  We also give a Bruhat interval characterization
of positroids.}

\keywords{positroids, decorated permutations, pattern avoidance,
Grassmannian}

\usepackage[backend=bibtex]{biblatex}
\begin{document}

\maketitle

\section{Introduction}
\label{sec: intro}

\textit{Positroids} are an important family of realizable matroids
originally defined by Postnikov in \cite{Postnikov.2006} in the
context of the totally nonnegative part of the Grassmannian
variety. These matroids and the totally positive part of the
Grassmannian variety have played a critical role in the theory of
cluster algebras and soliton solutions to the KP equations and have
connections to statistical physics, integrable systems, and scattering
amplitudes
\cite{AHBCGPT.2016,Rietsch.2006,williams2021positive}.
Positroids are closed under restriction, contraction, duality, and
cyclic shift of the ground set, and furthermore, they have particularly
elegant matroid polytopes \cite{Ardila-Rincon-Williams}.

\textit{Positroid varieties} were studied by Knutson, Lam, and Speyer
in \cite{KLS}, building on the work of Postnikov \cite{Postnikov.2006} and Rietsch \cite{Rietsch.2006}.  They are
homogeneous subvarieties of the complex Grassmannian variety $\Gkn$
that are defined by determinantal equations determined by the bases
of a positroid. They can also be described as projections of
Richardson varieties $X_{u} \cap X^{v}$ in the complete flag manifold
to $\Gkn$.  These varieties have beautiful geometric, representation
theory, and combinatorial connections \cite{KLS2,Paw}.  See the
background section for notation and further background.

The positroids $\M$ of rank $k$ on a ground set of size $n$ are in
bijection with many different combinatorial objects
\cite{Oh,Postnikov.2006}, including 
\begin{enumerate}
\item decorated permutations $\w$ on $n$ elements with $k$ anti-exceedances, 
\item Grassmann necklaces $(I_{1},\dotsc , I_{n}) \in 
\binom{[n]}{k}^{n}$,  and 
\item intervals $[u,v]$ in Bruhat order on $S_{n}$ such that $v$ is a 
$k$-Grassmannian permutation.  
\end{enumerate}
Here, a decorated permutation $\w$ on $n$ elements is a permutation $w \in
S_n$ together with an orientation clockwise or counterclockwise, denoted $\clockwise{i}$ or $\counterclockwise{i}$ respectively, on
the fixed points of $w$. In addition to these, there are bijections to juggling sequences, \reflectbox{L}-diagrams, equivalence classes of plabic graphs, and bounded affine permutations
\cite{Ardila-Rincon-Williams,KLS,Postnikov.2006}. In \Cref{sec:background}, we will sketch the relevant bijections and terminology.

Many of the properties of positroid varieties can be "read off" from one
or more of these bijectively equivalent definitions.  Thus, we will
index a positroid variety $\Pi_{\M}=\Pi_{\w}=\Pi_{[u,v]}$, depending on the relevant context.  For example, the codimension of $\Pi_{\M}$ in $\Gkn$ 
is easy to read off from the decorated permutation as follows.

Let $\SSnk$ be the set of decorated
permutations on $n$ elements with $k$ anti-exceedances. The \textit{chord diagram} $D(\w)$ of $\w \in
\SSnk$ is constructed by placing the numbers $1, 2,\ldots, n$ on $n$
vertices around a circle in clockwise order, and then, for each $i$,
draw a directed arc from $i$ to $w(i)$ with a minimal number of
crossings between distinct arcs while staying completely inside the
circle.  The arcs beginning at fixed points should be drawn clockwise
or counterclockwise according to their orientation in $\w$.

An \textit{alignment} in $D(\w)$ is a pair of directed edges $(i
\mapsto w(i), j \mapsto w(j)) $ which can be drawn as distinct
noncrossing arcs oriented in the same direction.  A pair of
directed edges $(i \mapsto w(i), j \mapsto w(j)) $ which can be drawn
as distinct noncrossing arcs oriented in opposite directions is
called a \textit{misalignment}.  A pair of directed edges which must
cross if both are drawn inside the cycle is called a  \textit{crossing}
\cite[Sect. 5]{Postnikov.2006}.  Let $\aligns(\w)$ denote the set of
alignments of $D(\w)$.

\begin{Ex} \label{Example:intro} Let $\w = 57 \counterclockwise{3}6492\clockwise{8}1$ be the decorated permutation with a
counterclockwise fixed point at $3$ and a clockwise fixed point at
$8$.  The chord diagram for $\w$ is the following.
\begin{center}
\includegraphics[height=4cm]{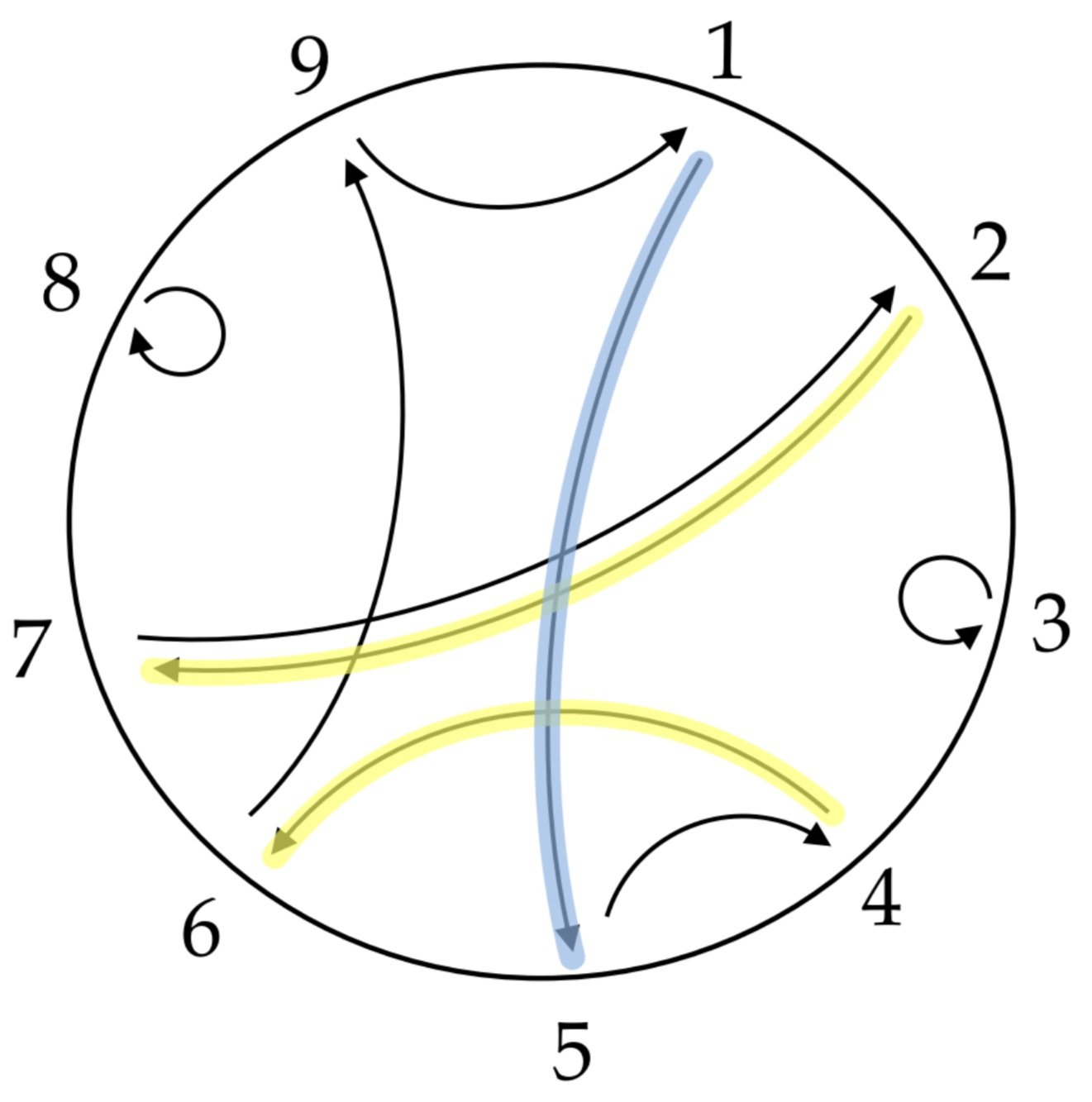}
\end{center}
Here, for example, $(2\mapsto 7, 4 \mapsto 6)$ highlighted in yellow is an alignment,
$(9\mapsto 1, 4 \mapsto 6)$ is a misalignment, and both $(1\mapsto 5,
2 \mapsto 7)$ and $(1\mapsto 5, 5 \mapsto 4)$ are crossings.  Note,
$(4\mapsto 6, 8 \mapsto 8)$ is an alignment, and $(4\mapsto 6, 3
\mapsto 3 )$ is a misalignment.  
\end{Ex}

\begin{Th}\label{thm:Alignments} \cite{KLS,Postnikov.2006} For any
decorated permutation $\w \in \SSnk$ and associated Bruhat interval
$[u,v]$, the codimension of $\Pi_{\w}$ in $\Gkn$ is
\begin{equation}\label{eq:codimPI}
\mathrm{codim}(\Pi_{\w}) =  \#\aligns(\w) = k(n-k) - [\ell(v) - \ell(u)].
\end{equation}
\end{Th}

Schubert varieties in the flag variety are indexed by permutations and
are closely related to positroid varieties, as explained in
\Cref{sec:background}.  Smoothness of Schubert varieties is
completely characterized by pattern avoidance of the corresponding
permutations \cite{Lak-San}. When studying the partially asymmetric
exclusion process and its surprising connection to the Grassmannian,
Sylvie Corteel posed the idea of considering patterns in decorated
permutations \cite{Corteel}.  This suggestion was the foundation for
the present work.

We use the explicit equations defining a positroid
variety in $\Gkn$ to determine if the variety is smooth or singular.
In general, a variety $X$ defined by polynomials $f_1, \ldots, f_s$ is
\textit{singular} if there exists a point $x \in X$ such that the
Jacobian matrix, $Jac$, of partial derivatives of the $f_i$ satisfies
$\text{rank}(Jac|_x) < \text{codim}\ X$. It is \textit{smooth} if no
such point exists.  The value $\text{rank}(Jac|_x)$ is the codimension
of the tangent space to $X$ at the point $x$.  Thus,
$\text{rank}(Jac|_x) < \text{codim} X$ implies the dimension of the
tangent space to the variety $X$ at $x$ is strictly larger than the
dimension of the variety $X$, hence $x$ is a singularity like a cusp on a
curve.  In the case of a positroid variety $\Pi_{\w}$,
\Cref{thm:Alignments} implies that a point $x \in \Pi_{\w}$ is a
singularity of $\Pi_{\w}$ if
\begin{align}
\text{rank}(Jac|_x) \; < \; \text{codim} \; \Pi_{\w} \; = \; \#\aligns(\w).
\end{align}

Our first main theorem reduces the problem of finding singular points
in a positroid variety to checking the rank of the Jacobian only at a
finite number of $T$-fixed points.  For any $J = \{j_1, \ldots, j_k\} \subseteq [n]$, let $A_{J}$ be the element in
$\Gkn$ spanned by the elementary row vectors $e_{i}$ with $i \in J$,
or equivalently the subspace represented by a $k\times n$ matrix with
a $1$ in cell $(i,j_i)$ for each $j_i \in J$ and zeros everywhere
else.  These are the $T$-fixed points of $\Gkn$, where $T \subset \text{GL}(n)$ is the set
of invertible diagonal matrices over $\mathbb{C}$. The reduction follows from the
decomposition of $\Pi_{[u,v]}$ as a projected Richardson variety. Every point $A \in \Pi_{[u,v]}$ lies in the projection of some
intersection of a Schubert cell with an opposite Schubert variety
$C_{y}\cap X^{v}$ for $y\in [u,v]$. In particular,
if  $y=y_{1} y_{2}\cdots y_{n} \in [u,v]$ in one-line notation and we define
$y[k] :=\{y_{1},y_{2},\dots, y_{k}\}$, then $A_{y[k]}$ is in the
projection of $C_y \cap X^v$.

\begin{Th}\label{thm:bounded.below}
Assume $A \in \Pi_{[u,v]}$ is the image of a point in $C_{y} \cap X^v$
projected to $\Gkn$ for some $y \in [u,v]$. Then the codimension of the tangent space to
$\Pi_{[u,v]}$ at $A$ is bounded below by $\rank(Jac|_{A_{y[k]}})$.
\end{Th}

\Cref{thm:bounded.below} indicates that the $T$-fixed points of the
form $A_{y[k]}$ such that $y \in [u,v]$ are key for understanding the singularities
of $\Pi_{[u,v]}$.  In fact, the equations determining $\Pi_{[u,v]}$
and the bases of the positroid $\M$ associated with the interval $[u,v]$
can be determined from the permutations in the interval by the following
theorem.  Our proof of the following theorem depends on \Cref{Lemma:
BS}. It also follows from \cite[Lemma 3.11]{KW2015}.

\begin{Th}  \label{Th: positroid = initial sets} \normalfont Let
$\w \in \SSnk$ have associated Bruhat interval $[u,v]$ and
positroid $\M$. Then $\M$ is exactly the collection of initial sets of permutations in the Bruhat interval $[u,v]$,
\[
\M = \{y[k] \; : \; y \in [u,v]\}.
\]
\end{Th}

Our next theorem provides a method to compute the rank of the
Jacobian of $\Pi_{{[u,v]}}$ explicitly at the $T$-fixed points.  Therefore, we can also compute the
dimension of the tangent space of a positroid variety at those points.

\begin{Th} \label{thm:tangent.space.dim} Let
$\w \in \SSnk$ have associated Bruhat interval $[u,v]$ and
positroid $\M$. For any  $y\in [u,v]$, the codimension of the tangent space to $\Pi_{[u, v]}
\subseteq \Gkn$ at
$A_{y[k]}$ is
\begin{equation}\label{eq:tangent.space.codim}
\rank (Jac|_{A_{y[k]}}) \; = \; \# \Big\{I \in \binom{[n]}{k} \setminus \M \; : \; |I \cap y[k]| = k-1
\Big\}.
\end{equation}
\end{Th}

The formula in \eqref{eq:tangent.space.codim} is reminiscent of the
\textit{Johnson graph} $J(k,n)$ with vertices given by the $k$-subsets
of $[n]$ such that two $k$-subsets $I,J$ are connected by an edge precisely
if $|I\cap J|=k-1$.  For a positroid $\M \subseteq \binom{[n]}{k}$, let $J(\M)$ denote the induced subgraph of the
Johnson graph on the vertices in $\M$.  We call $J(\M)$ the
\textit{Johnson graph of} $\M$.  Note, the Johnson graph of $\M$ is closely
related to the Basis Exchange Property for matroids.
\Cref{thm:tangent.space.dim} implies $J(\M)$ encodes aspects of the
geometry of the positroid varieties like the Bruhat graph in the
theory of Schubert varieties \cite{carrell94}. 

To state our main theorem characterizing smoothness of positroid
varieties, we need to define two types of patterns that may occur in a
chord diagram.  Given an alignment $(i \mapsto w(i), j \mapsto w(j))
$ in $D(\w)$, if there exists a third arc $(h \mapsto w(h))$ which
forms a crossing with both $(i \mapsto w(i))$ and $(j \mapsto w(j))$,
we say $(i \mapsto w(i), j \mapsto w(j)) $ is a \textit{crossed
alignment} of $\w$. In the example above, $(2\mapsto 7, 4 \mapsto 6)$
is a crossed alignment; this alignment is crossed for instance by $(1
\mapsto 5)$, highlighted in blue. We call a permutation of the
form $w(i)=i+t \, (mod\ n)$ for some fixed integer $1\leq t \leq n$ a
\textit{spirograph permutation}, and we call its chord diagram a
\textit{spirograph}. We think of alignments, crossings, crossed
alignments, and spirographs as subgraph patterns for decorated permutations.

%

\begin{Th} \label{The: Main theorem} For $\w \in \SSnk$
with associated positroid  $\M$,  the following are
equivalent.
\begin{enumerate}
\item The positroid variety $\Pi_\M$ is smooth.
\item For every $J \in \M$,  $\# \{I \in \M \; : \; |I \cap J| = k-1 \} \; = \; k(n-k) - \# \aligns(\w).$
\item The graph $J(\M)$ is regular, and each vertex
has degree  $k(n-k)-\# \aligns(\w).$
\item The decorated permutation $\w$ has no crossed alignments.
\item The chord diagram of $\w$ is a union of spirographs on a
noncrossing partition of $[n]$.
\end{enumerate}
\end{Th}

In this extended abstract, we outline the steps for the proof of
\Cref{The: Main theorem} and leave the details of the remaining proofs as well as enumerative results to the
forthcoming paper \cite{Billey-Weaver} and OEIS entries A349413, A349456, A349457, and A349458.  In
\Cref{sec:background}, we provide background material and define
our notation.  In \Cref{Sec:
reductions}, we reduce the proof of \Cref{The: Main theorem} to
derangements.   We provide further reductions which state that $\Pi_{\w}$ is smooth if
and only if the positroid varieties defined by flipping, inverting, or
rotating $D(\w)$ are also smooth.  Finally, we
use the fact that we can flip and rotate a chord diagram with a
crossed alignment so that the crossing arc has its tail at 1 and
crosses the alignment from starboard to port when the alignment is
viewed as a boat moving forward, as shown in the highlighted crossed
alignment in \Cref{Example:intro}.  In this configuration, we show
that the $T$-fixed point associated with the anti-exceedance set of
$\w$ is singular in $\Pi_{\w}$.  If no crossed alignment exists, we
show that $D(\w)$ is a disjoint union of spirographs, which leads to a
regular Johnson graph of $\M$.

\section{Background}\label{sec:background}

We begin by giving notation and some background on several
combinatorial objects and theorems from the literature. These objects
will be used to index the varieties discussed throughout the paper. We
will then introduce our notation for the flag variety, the
Grassmannian, Schubert varieties, Richardson varieties, and positroid
varieties.  See Fulton's book on Young Tableaux for further background
\cite{Ful}.

\subsection{Combinatorial objects}

For integers $i \leq j$, let $[i,j]$ denote the set
$\{i,i+1,\ldots,j\}$, and write $[n] := [1,n]$. Let $\binom{[n]}{k}$
be the set of $k$-element subsets of $[n]$. Call $J \in
\binom{[n]}{k}$ a $k$-subset of $[n]$. Let $A_J$ be the
$k \times n$ matrix whose restriction to the columns indexed by $J$ is
the identity matrix $I_k$, and whose other entries are zeros.

\begin{Ex} For $J = \{2,4,8\} \in \binom{[9]}{3}$ , $A_J$ is the matrix
\[
A_J = 
\begin{bmatrix}
0 & 1 & 0 & 0 & 0 & 0 & 0 & 0 & 0\\
0 & 0 & 0 & 1 & 0 & 0 & 0 & 0 & 0\\
0 & 0 & 0 & 0 & 0 & 0 & 0 & 1 & 0
\end{bmatrix}.
\]
\end{Ex}

\medskip

Define the \textit{Gale partial order}, $\preceq$, on $k$-subsets of
$[n]$ as follows. Let $I = \{i_1 < \cdots < i_k\}$ and $J = \{j_1 <
\cdots < j_k\}$. Then $I \preceq J$ if and only if $i_h \leq j_h$ for
all $h \in [k]$. This partial order is known by many other names; we
are following \cite{Ardila-Rincon-Williams} for consistency.

A \textit{matroid} of rank $k$ on $[n]$, defined by its bases, is a
nonempty subset $\mathcal{M} \subseteq \binom{[n]}{k}$ satisfying the
following Basis Exchange Property: if $I, J \in \mathcal{M}$ such that
$I \neq J$ and $a \in I \setminus J$, then there exists some $b \in J
\setminus I$ such that $(I \setminus \{a\}) \cup \{b\} \in
\mathcal{M}$. Compare the notion of matroid basis exchange to basis
exchange in linear algebra.

For example, let $A$ be a full rank $k \times n$ matrix. The \textit{matroid of $A$} is the set
\[
\mathcal{M}_A := \Big\{J \in \binom{[n]}{k} \; : \; \Delta_J(A) \neq 0\Big\},
\]
where $\Delta_J(A)$ is the determinant of the $k \times k$ submatrix of $A$ lying in column set $J$. 

\begin{Ex}\label{ex:matroid} The matroid of
\[
A = \begin{bmatrix}
0 & 3 & 1 & -2 &  2 & 0 \\
0 & 0 & 0 &  1 & -1 & 1
\end{bmatrix} 
\]
is \[
\mathcal{M}_A = \{ \{2,4\}, \{2,5\}, \{2,6\}, \{3,4\}, \{3,5\}, \{3,6\}, \{4,6\}, \{5,6\}\} \subseteq \binom{[6]}{2}.
\]
\end{Ex}

Let $S_n$ be
the set of permutations of $[n]$, where we think of a permutation as a
bijection from the set $[n]$ to itself. For $w \in S_n$, let
$w_i=w(i)$, and write $w$ in \textit{one-line notation} as
$w=w_1w_2\cdots w_n$.  A permutation with no fixed points $i=w(i)$ is
a \textit{derangement}. The
permutation matrix $M_w$ of $w$ is the $n \times n$ matrix that has a
1 in cell $(i,w_i)$ for each $i \in [n]$ and zeros elsewhere. The
\textit{length of} $w \in S_n$ is the number of inversions in $w$,
\[
\ell(w) := \# \{(i,j) \; : \; i < j \; \text{and} \; w(i) > w(j)\}.
\]

\begin{Ex} For $w = 3124$, the length of $w$ is $\ell(w) = 2$, and $M_w$ is the matrix
\[
M_{3124} = 
\begin{bmatrix}
0 & 0 & 1 & 0 \\
1 & 0 & 0 & 0 \\
0 & 1 & 0 & 0 \\
0 & 0 & 0 & 1
\end{bmatrix}.
\]
\end{Ex}

\noindent Note, the permutation matrix of $w^{-1}$ is $M_{w}^{T}$.
Permutation multiplication is given by function
composition so that if $wv=u$, then $w(v(i))= u(i) $. Hence, $M_w^TM_v^T= M_u^T$.

For $1\leq k \leq n$, write $S_k \times S_{n-k}$ for the subgroup of
$S_n$ consisting of permutations that send $[k]$ to $[k]$ and
$[k+1,n]$ to $[k+1,n]$. For $0 \leq k \leq n$, a permutation $w \in S_n$ is a
$k$-\textit{Grassmannian permutation} if $w_1 < \cdots < w_k$ and
$w_{k+1} < \cdots < w_n$. This is equivalent to saying
that $w$ is the minimal length element of its coset $w \cdot (S_k
\times S_{n-k})$. The permutation $w=35124$ is 2-Grassmannian.

\begin{Def}\label{th:tableau.criterion} \cite[Chapter 2]{BjornerBrenti}
For $u,v \in S_{n}$,\  $u\leq v$ in \textit{Bruhat order} if
$u[i]\preceq v[i]$ for all $i \in [n]$.   
\end{Def}

For each $u\leq v$ in Bruhat order, the \textit{interval} $[u,v]$ is
defined to be $[u,v] :=\{y\in S_{n}\ : \ u \leq y\leq v \}$.  The
intervals $[u,v]$ where $v$ is a $k$-Grassmannian permutation are key
to this work.  In this case, the following simpler criterion for Bruhat
order follows from work of Bergeron--Sottile \cite[Theorem
A]{BergSott}.

\begin{Th}\label{Lemma: BS}  Let $u,v \in
S_n$, and assume $v$ is $k$-Grassmannian.  Then $u \leq v$ if and only if
\begin{enumerate}
\item[(i)] for every $1 \leq j \leq k$, we have $u(j) \leq v(j)$, and
\item[(ii)] for every $k < m \leq n$, we have  $u(m) \geq v(m)$.  
\end{enumerate}
\end{Th}

\subsection{Grassmannian, Flag, and Richardson Varieties} \label{Sec: Fl(n) and Gr(k,n)}

For $0 \leq k \leq n$, the points in the \textit{Grassmannian variety},
$\Gkn$, are the $k$-dimensional subspaces of $\C^n$. Up to left
multiplication by a matrix in $\GL_k$, we may represent $V \in \Gkn$
by a full rank $k \times n$ matrix $A_V$ such that $V$ is the row span
of $A_V$. Let $\mathrm{Mat}_{kn}$ be the set of full rank $k \times n$
matrices. One may think of $\Gkn$ as the cosets $\GL_k \backslash
\mathrm{Mat}_{kn}$.  The Grassmannian varieties are smooth manifolds. This includes the case when $k=n=0$, in which case $\Gkn$ is the single point $(0)$.

Let $\mathcal{F}\ell(n)$ be the \textit{complete flag variety} of
nested subspaces of $\C^n$.  A complete flag $V_{\bullet} =(V_{1},
\ldots , V_{n})$ can be represented as an
invertible $n \times n$ matrix where the row span of the first $j$
rows corresponds with the $j^{th}$ subspace in the flag.  For a subset $J
\subseteq [n]$, let Proj$_J : \C^n \rightarrow \C^{|J|}$ be the
projection map onto the indices specified by $J$. Then for every
permutation, $w \in S_n$, there is a Schubert cell $C_w$ and an
opposite Schubert cell $C^w$ in $\mathcal{F}\ell(n)$ defined by
\begin{align*}
&C_w = \{V_\bullet \in \mathcal{F}\ell(n) \; : \; \dim(\text{Proj}_{[j]}(V_i)) = |w[i] \cap [j]| \;\text{for all $i,j$}\} \text{ and} \\
&C^w = \{V_\bullet \in \mathcal{F}\ell(n) \; : \; \dim(\text{Proj}_{[n-j+1,n]}(V_i)) = |w[i] \cap [n-j+1,n]| \;\text{for all $i,j$}\}.
\end{align*}
The \textit{Schubert variety} $X_w$ is the closure of $C_w$ in the
Zariski topology on $\mathcal{F}\ell(n)$, and similarly, the
\textit{opposite Schubert variety} $X^w$ is the closure of $C^w$.
Bruhat order determines which Schubert cells are in a Schubert variety, 
\begin{equation}\label{eq:schubs}
 X_w = \bigcupdot_{y \geq w} C_{y} \ \text{ and } \ X^{w} = \bigcupdot_{v
\leq w} C^{v}.
\end{equation}

For permutations $u$ and $v$ in $S_n$, with $u \leq v$, the
\textit{Richardson variety} is a nonempty variety in
$\mathcal{F}\ell(n)$ and is defined as the intersection $X_u^v := X_u
\cap X^v$. Then $\dim X_u^v = \ell(v) - \ell(u)$. 
The decompositions of $X_u$ and $X^v$ into Schubert cells and opposite Schubert cells yield
\[
X_u^v = {\displaystyle \bigcupdot_{u \leq y \leq v} (C_y \cap X^v)}
=
\Big({\displaystyle \bigcupdot_{y \geq u}} C_y \Big) \cap
\Big({\displaystyle \bigcupdot_{t \leq v}} C^t\Big).
\]

\subsection{Positroids and Positroid Varieties}\label{sub:positroids}

Postnikov and Rietsch considered an important cell decomposition of
the totally nonnegative Grassmannian
\cite{Postnikov.2006,Rietsch.2006}.  The term \textit{positroid} does
not appear in either paper, but has become the name for the matroids that index the nonempty
matroid strata in that cell decomposition.  They also individually considered the
closures of those cells, which determines an analog of Bruhat order.
The cohomology classes for these cell closures was investigated by
Knutson, Lam, and Speyer \cite{KLS,KLS2} and Pawlowski \cite{Paw}.

\begin{Def} \label{Def: positroid}
\normalfont Let $\Gkn^{tnn}$ be the points in $\Gkn$ that can each be represented by
a real valued $k\times n$ matrix $A$ such that every minor $\Delta_I(A)$ satisfies
$\Delta_{I}(A)\geq 0$.  A \textit{positroid} is a matroid of the form
$\M_A$ for some matrix $A \in \Gkn^{tnn}$.
\end{Def}

Postnikov also made the following definitions in \cite[Sect
16]{Postnikov.2006}.  Given a decorated permutation $\w$, as
defined in the introduction, call $i \in [n]$ an
\textit{anti-exceedance} of $\w$ if $i<w^{-1}(i)$ or if $w(i) = \clockwise{i}$ is
a clockwise fixed point. Fix $r \in [n]$. Let $<_r$ be the
shifted linear order on $[n]$ given by
$r <_r (r+1) <_r \cdots <_r n <_r 1 <_r \cdots <_r (r-1).$
The \textit{shifted anti-exceedance set $I_r(\w)$} of $\w$ is the
anti-exceedance set of $\w$ with respect to the shifted linear order
$<_r$ on $[n]$,
\[
I_r(\w) = \{i \in [n] \; : \; i <_r w^{-1}(i)  \; \text{or} \; w(i) = \clockwise{i} \}.
\]
Thus, $I_1(\w)$ is the set of anti-exceedances of $\w$. Recall that $\SSnk$
is the set of decorated permutations with anti-exceedance set
$I_1(\w)$ of size $k$. The \textit{Grassmann necklace} associated with
$\w$ is $(I_1(\w), \ldots, I_n(\w))$.  By construction,
$|I_1(\w)|=\cdots=|I_n(\w)|=k$.

Using the shifted linear order on $[n]$, we may define the
\textit{shifted Gale order} $\prec_r$ on
$\binom{[n]}{k}$. Specifically, if $I, J \in \binom{[n]}{k}$, where $I = \{i_1 <_r \cdots <_r i_k\}$ and $J = \{j_1 <_r
\cdots <_r j_k\}$, then $I \preceq_r J$ if $i_h \leq_r j_h$ for all $h
\in [k]$. The positroid associated with a decorated permutation can be defined using the shifted anti-exceedance sets and the shifted Gale orders.

\begin{Th}\label{th:decoratedpermtopositroid}\cite{Oh,Postnikov.2006}
For $\w \in \SSnk$, the set
\begin{align} \label{Def: positroid def from anti-ex sets}
\M(\w) \; := \; \Big\{I \in \binom{[n]}{k} \; : \; I_r(\w) \preceq_r I \text{ for all } r \in [n]\Big\}
\end{align}
is a positroid. Conversely, for every positroid $\M$ of rank $k$ on ground set $[n]$,
there exists a unique decorated permutation $\w\in \SSnk$ such that the 
sequence of minimal elements in the shifted Gale order on the subsets
in $\M$ is the Grassmann necklace of $\w$.  
\end{Th}

To find the interval $[u,v]$ corresponding to $\w$ following
\cite{KLS}, compute the set $I_{1} = I_1(\w)$, say
$|I_{1}|=k$. Let $w$ be the permutation associated with $\w$ by
forgetting the decorations. Let $v$ be the unique $k$-Grassmannian
permutation such that $\{v_{1},\ldots ,v_{k} \} =w^{-1}(I_{1})$. Then 
$wv=u$, so $I_{1}=u[k]=\{u_{1},\ldots , u_{k} \}$.  To recover $\w$
from $[u,v]$, compute $w=uv^{-1}$ and orient all of the fixed points
of $w$ which are in $u[k]$ to be clockwise, and orient all others to be
counterclockwise.

For example, for the decorated permutation $\w = 57 \counterclockwise{3}6492\clockwise{8}1$ in \Cref{Example:intro}, we
have $k=4$, the Grassmann necklace is
\[
(I_{1},\ldots, I_{9})= (\{ 1248 \}, \{ 2458 \}, \{ 4578 \}, \{ 4578\} ,
\{ 5678 \}, \{ 4678 \} , \{ 4789 \}, \{ 2489 \} , \{ 2489 \}), 
\]
and the corresponding interval $[u,v]$ has $u=428157369$ and
$v=578912346$.  The corresponding positroid has 22 elements.

In \Cref{ex:matroid}, the matrix $A$ has all nonnegative $2\times 2$
minors, so the associated matroid is a positroid.  The minimal elements
in shifted Gale order are $(\{2 4\}, \{2 4\},$ $\{3 4\},$ $\{4 6\}, \{5
6\},$ $\{2 6\})$, which is the Grassman necklace for the decorated
permutation $\counterclockwise{1}36524$.  The associated Bruhat interval is
$[241365, 561234]$.

As mentioned in the introduction, there are many other objects in
bijection with positroids and decorated permutations.  We refer the
reader to \cite{Ardila-Rincon-Williams} for a nice survey of many
other explicit bijections.

Let $\pi_k : \mathcal{F}\ell(n) \rightarrow Gr(k,n)$ be the projection
map which sends a flag $V_\bullet = (0 \subset V_1 \subset \cdots
\subset V_n)$ to the $k$-dimensional subspace $V_k$. Identifying a
full rank $n \times n$ matrix $M$ with the point it represents in
$\mathcal{F}\ell(n)$, then $\pi_k(M)$ denotes the span of the top $k$
rows of $M$. 

\begin{ThmDef} \label{Def: positroid as projection}\cite[Thm 5.1]{KLS}
Given a decorated permutation $\w \in \SSnk$ along with its
associated Bruhat interval $[u,v]$ and positroid $\M\subseteq
\binom{[n]}{k}$, the following are equivalent definitions of the
\textit{positroid variety} $\Pi_{\w} = \Pi_{[u,v]} = \Pi_{\M}$.  
\begin{enumerate}
\item The positroid variety $\Pi_{\M}$ is the homogeneous subvariety of
$Gr(k,n)$ whose vanishing ideal is generated by the Pl{\"u}cker
coordinates $\{\Delta_I \; : \; I \notin \M\}$.
\item The positroid variety $\Pi_{[u,v]}$ is the projection of the
Richardson variety $X_u^v$ to $Gr(k,n)$, so $\Pi_{[u,v]}=\pi_k(X_u^v)$.
\end{enumerate}
\end{ThmDef}

\section{Outline of Proofs}\label{Sec: reductions}

In this section, we outline the steps for reducing the proof of \Cref{The: Main theorem} to equivalence classes of derangements
under flip, inverses, and rotation. \Cref{Th: positroid = initial sets}
leads to the first reduction, allowing us to effectively ignore
fixed points. In particular, $\w$ has a fixed point $j$ if and only if
the corresponding $[u,v]$ has an index $i$ such that $u(i)=v(i)=j$.
In this case, every $y \in [u,v] $ has $y(i)=j$, so the
effect of removing the fixed point $j$ from $\w$ is easily described
in terms of the Bruhat interval, and hence the
positroid.  In turn, the equations defining the positroid variety
for $\w$ and the positroid variety indexed by the decorated permutation $\vv$ obtained by removing $j$ from $\w$ are
closely related.   

\begin{Lemma} \label{Lemma: removing fixed points}
Let $\w \in \SSn$ be a decorated permutation with fixed point $i$. Let $\vv \in \SSnminus$ be obtained from $\w$ by deleting the loop at $i$ from $D(\w)$. Then $\Pi_{\w}$ is smooth if and only if $\Pi_{\vv}$ is smooth.
\end{Lemma}

If $\w$ fixes every point in $[n]$, then $\Pi_{\w}$ is a single point, and hence is smooth. Otherwise, by \Cref{Lemma: removing fixed points}, instead of considering $\w \in \SSn$, we may consider the
derangement obtained by deleting all the fixed points in $\w$. Thus, we may restrict our attention to derangements. For the remainder, we give all
results for derangements in $S_n$, and we drop the decorated
permutation notation. In \cite{Billey-Weaver}, we further reduce to \textit{stabilized interval-free permutations}, found in \cite[A075834]{OEIS}.

\begin{Lemma} \label{Lemma: invariances} Let $w$ be a
derangement in $S_n$, and let $w'$ be a derangement obtained from
$w$ by (1) rotating the chord diagram of $w$, (2) reflecting $D(w)$ across the vertical axis, or (3) reversing the direction of all arcs. Then, in any case, $\Pi_w$ is smooth if and only if $\Pi_{w'}$ is smooth.
\end{Lemma}

The chord diagram of $w \in S_n$ is a disjoint union of $m$
graphs if the arcs can be partitioned into $m$ parts
such that no arcs from distinct parts form a crossing.  In this case,
one may consider the restriction of the diagram to any of these parts
and the corresponding decorated permutations.  

\begin{Lemma} \label{Lemma: connected components}
Suppose $w \in S_n$ is a derangement whose chord diagram can be
partitioned into a disjoint union of $m$ graphs. Let $w^{(1)}, \ldots,
w^{(m)}$ be derangements corresponding to these parts. Then $\Pi_w$ is
smooth if and only if all of the $\Pi_{w^{(i)}}$ are smooth.
\end{Lemma}

The proofs of the lemmas above rely on \Cref{thm:tangent.space.dim}
and the Johnson graphs of the corresponding positroids.  Explicit maps
on edges in the Johnson graphs are given for each type of reduction.

Next, we give the outline of the proof of \Cref{The:
Main theorem}.  The first three items are equivalent from the
definition of smooth, \Cref{thm:bounded.below}, and \Cref{thm:tangent.space.dim}.  Next, we show  parts 
 $ 4 \iff 5,\  5\Rightarrow 1$, and $
2 \Rightarrow 4$.

A spirograph permutation $w$ of the form $w(i) = i + t \,
(mod \, n)$ has no alignments, so it has no crossed
alignments.  For the reverse direction, one may reduce to a
single connected component of $D(w)$ as in \Cref{Lemma: connected components}
and assume that $w$ does not have the form $w(i) = i+t \, (mod \,
n)$. Then, one can directly find a crossed alignment.  Thus, $ 4
\iff 5$.

The fact that $w$ satisfying $w(i) = i + t \, (mod \, n)$ has no
alignments also implies that $\Pi_w = Gr(k,n)$ for the appropriate
value of $k$. Thus, $\Pi_w$ is smooth. \Cref{Lemma: connected
components} then provides the implication $5 \Rightarrow 1$ in
\Cref{The: Main theorem}.

The final step in proving \Cref{The: Main theorem} is achieved by
proving $2 \Rightarrow 4$ by the contrapositive. In particular, given
a derangement $w$ with a crossed alignment, we identify a set $J \in
\M$ such that the number of non-bases $I \in \QQ = \binom{[n]}{k} \setminus \M$
for which $|I \cap J| = k-1$ is strictly less than the number of
alignments of $w$.  By applying \Cref{Lemma: invariances} as
necessary, we may assume that the crossing arc of the crossed
alignment passes from starboard to port and has its tail at 1, e.g., the highlighted crossed alignment in \Cref{Example:intro}. In this
regime, consider the anti-exceedance set $J = I_1(w)$.

\begin{Lemma} \normalfont \label{Lemma: anti-exchange pair condition} For a derangement $w \in S_n$ with corresponding positroid $\M$, let $J = I_1(w)$, and suppose that $a \in J$ and $b \notin J$. Then $I = (J \setminus \{a\}) \cup \{b\}$ is in $\M$ if and only if $a < b$ and for every $r \in [a+1,b]$, the following two conditions hold.
\begin{enumerate}
\item There exists some $x \in [a,r-1]$ such that $w^{-1}(x) \geq r$.
\item There exists $y \in [r,b]$ such that $w^{-1}(y) \leq r-1$.
\end{enumerate}
\end{Lemma}

\Cref{Lemma: anti-exchange pair condition} gives an exact condition on the sets $I \in \QQ$ for which $|I \cap J| = k-1$. Any such set $I$ corresponds to the pair $(a,b)$ such that $I = (J \backslash \{a\}) \cup \{b\}$. We call such a  pair $(a,b)$ an \textit{anti-exchange pair} for $J$ and define a map $\Psi_w$ from the anti-exchange pairs for $J$ to $Alignments(w)$. In the case that $(a,b)$ is an anti-exchange pair for $J$ with $b < a$, then $\Psi_w(a,b) = (w^{-1}(b) \mapsto b, w^{-1}(a) \mapsto a)$. Otherwise, $a < b$ and at least one of condition 1 or 2 in \Cref{Lemma: anti-exchange pair condition} is not satisfied for some $r \in [a+1,b]$.

Conditions 1 and 2 of the lemma are symmetric. The map $\Psi_w$
inherits this symmetry, and thus we will only give the details of the
map for anti-exchange pairs that do not satisfy condition (1). In this
case, by \Cref{Lemma: anti-exchange pair condition}, there exists some $r
\in [a+1,b]$ such that for all $x \in [a,r-1]$, $w^{-1}(x) \leq
r-1$. Choose r to be minimal such that condition 1 is not satisfied. By starting at $b$ and tracing in reverse the cycle containing $b$, there will be some first element $c \neq b$ such that $c \geq r$. Then $w(c) \in [c+1,n] \cup [a-1]$ and $w^{-1}(a) \in [a+1,r-1]$, since $a$ is an anti-exceedance. Therefore, we choose $\Psi_w(a,b) = (c \mapsto w(c), w^{-1}(a) \mapsto a)$. 

\begin{Lemma} \label{Lemma: injective map} Let $w \in S_n$ be a derangement. The map $\Psi_w$ from anti-exchange pairs of $I_1(w)$ to $Alignments(w)$ is injective. Furthermore, if $w$ has a crossed alignment with crossing arc passing from starboard to port and whose tail is at 1, then $\Psi_w$ is not surjective.
\end{Lemma}

In particular, in the latter case, there is no anti-exchange pair for
$I_1(w)$ that is mapped to the crossed alignment. Thus, if $w$ has a
crossed alignment, then $\# Alignments(w)$ is strictly larger than the
number of anti-exchange pairs for $I_1(w)$. Since $I_1(w) = u[k]$, where
$[u,v]$ is the Bruhat interval corresponding to $w$, then $I_1(w)
\in \M=\{y[k]: \ y \in [u,v] \}$ by \Cref{Th: positroid = initial sets}. Therefore, $I_1(w)$ is an element of $\M$ for which part 2 of  \Cref{The: Main theorem} is not satisfied. It then follows that $2 \Rightarrow 4$ in \Cref{The: Main theorem}.

\section{Acknowledgments}\label{Sec:Future}

We thank Herman Chau, Lauren Williams, Brendan Pawlowski, Stark Ryan, Joshua Swanson, and Allen Knutson for many insightful comments on this project.

\printbibliography


\end{document}